\documentclass[12pt]{amsart}
\theoremstyle{definition}

\theoremstyle{remark}

\numberwithin{equation}{section}

\begin{document}
\setlength{\baselineskip}{16pt}
\title{Integer Solutions of a Sequence of Decomposable Form Inequalities}
\author{K. Gy\H ory and Min Ru}
\address{Institute of Mathematics and Informatics, Lajos Kossuth University, 
H-4010 Debrecen, Pf 12, Hungary;
\newline
 Department of Mathematics, University of 
Houston, Houston, TX 77204}
\email{gyory@math.klte.hu, minru@math.uh.edu}
\thanks{
The first author was supported in part by the Hungarian 
Academy of Sciences  and by Grants 16975 and 25157 from the Hungarian 
National Foundation for Scientific Research. The second author was supported 
in part by NSF grant DMS-9596181}
\maketitle 

\section{Introduction}
Let $F({\bf X}) = F(X_0, \dots, X_m) \in {\bf Z}[{\bf X}]$ be a 
decomposable form, i.e. a homogeneous polynomial which 
factorizes into linear forms over
 ${\bar {\bf Q}}$. Assume that $q = \deg F > 2m$, and consider the 
decomposable form inequality 
\begin{equation}
\begin{split}
0 < |F({\bf x})| <  |{\bf x}|^{\nu} ~\mbox{in} ~{\bf x} 
= (x_{0}, \dots, x_{m}) \in 
{\bf Z}^{m+1},
\end{split}
\end{equation}
where $|{\bf x}| = \max_{0 \leq i \leq m} |x_i|$ and $\nu < q-2m$. For $m=1$, 
it follows from Roth's approximation theorem (cf. [Sch 2], p. 120) that if the linear factors of $F$ are pairwise non-proportional, then (1.1) has only
 finitely many solutions. Using his subspace theorem, W.M. Schmidt 
([Sch 1], [Sch 2]) generalized this for arbitrary $m$, under the assumptions 
that (i) any $m+1$ of the linear factors of $F$ are linearly independent
 over ${\bar {\bf Q}}$, and
that (ii) $F$ is not divisible by a form with rational coefficients of
 degree less than $m+1$. Later H.P. Schlickewei [Schl] extended this theorem 
to the case when the ground ring is an arbitrary finitely generated subring 
of ${\bf Q}$. These results have obvious applications to decomposable 
form equations of the form 
\begin{equation}
\begin{split}
F({\bf x}) = G({\bf x}) ~\text{in}~{\bf x}  \in  {\bf Z}^{m+1},
\end{split}
\end{equation}
where $G \in {\bf Z}[{\bf X}]$ is a non-zero polynomial of degree $< q-2m$ 
(cf. [Sch 2]). In the important special case when $G$ is a constant, 
the first author [Gy 1] proved the finiteness of the number of solutions of 
(1.2) under
the assumption (i) only, which is necessary in general. For the case
when  $G$ is a constant, 
there are also more general finiteness results, 
see 
 [Sch 2], [EGy 2], [EGy 4] and the references given there. 
The above mentioned results have further applications to resultant 
inequalities and resultant equations (cf. [Sch 1], [Sch 2], 
[Schl], [Gy 1], [Gy 3]). 

\bigskip
The purpose of this paper is to improve and generalize the above results of 
 [Sch 1], [Sch 2], [Schl] for a sequence of decomposable form inequalities
 over number fields (cf. Theorem 2). For a single decomposable form 
inequality, Theorem 2 implies the finiteness of the number of solutions of 
(1.1) without 
assuming (ii) (cf. Theorem 1). Our Theorem 1 and Theorem 2 do not remain 
valid in general for $\nu=q-2m$. As a consequence of Theorem 2,
finiteness results 
are established for decomposable form equations of 
 the form (1.2) over number fields (cf. Theorems 3, 4). 
Some applications
are also given to resultant inequalities (Theorems 5, 6).
Finally, we generalize our results for the case where $q>2m-l+1$ with 
$1 \leq l \leq m+1$, by showing that in this situation the set of solutions is 
contained in a finite union of subspaces of dimension at most $l$ 
(cf. Theorems 7, 8, 9).

\section{Notation and statement of the main results}
Let {\bf K} be an algebraic number field. Denote by {\bf M(K)} the set
of places of {\bf K} and write ${\bf M_{\infty}(K)}$ for the set of 
archimedean
places of {\bf K}.  For $\upsilon \in {\bf M(K)}$ denote 
by $|~~|_{\upsilon}$ the
associated absolute value, normalized such that 
$|~~|_{\upsilon}
= |~~|$ (standard absolute value) on {\bf Q} if $\upsilon$ is archimedean, whereas for $\upsilon$
non-archimedean $|p|_{\upsilon} = p^{-1}$ if $\upsilon$ lies above the 
rational 
prime p. Denote by ${\bf K}_{\upsilon}$ 
the completion of {\bf K} with respect to $\upsilon$ and by $d_{\upsilon} = [{\bf K}_{\upsilon}
: {\bf Q}_{\upsilon}]$ the local degree. We put $||~~||_{\upsilon} = |~~|_{\upsilon}^{d_{\upsilon}}.
$

\bigskip
For ${\bf x} = (x_0, \dots, x_n) \in {\bf K}^{m+1}$, we put $\|{\bf x}\|_{\upsilon} = \max _{0 \leq i \leq m}\|x_i\|_{\upsilon}$ and we denote by 
\[ h({\bf x}) = {1\over [{\bf K}:{\bf Q}]} \sum_{\upsilon \in {\bf M}({\bf K})} \log \|{\bf x}\|_{\upsilon} \]
the absolute logarithmic height of ${\bf x}$. Given a polynomial $P$ with
 coefficient in ${\bf K}$, we define $\|P\|_{\upsilon}$ and $h(P)$ as the 
$\|~\|_{\upsilon}$-value and absolute logarithmic height, respectively, of the
point whose coordinates are the coefficients of $P$. As is known, 
$h({\bf x})$ and $h(P)$ are independent of the choice of the 
field ${\bf K}$. Further, 
$h(\lambda {\bf x}) = h({\bf x})$ and $h(\lambda P) = h(P)$ for all 
$\lambda \in {\bar {\bf Q}}^*$.

\bigskip
Let S be a finite subset of {\bf M(K)} containing ${\bf M_{\infty}(K)}$. 
An element $x \in {\bf K}$ is said to be S-integer if
 $||x||_{\upsilon} \leq 1 $
for each $\upsilon \in M(K)-S$. Denote by $O_S$ the set of S-integers.
The units of $O_S$ are called $S$-units. They form a multiplicative group 
which is denoted by $O^*_S$.
For $\alpha \in {\bf K} - \{0\}$, let 
$N_S(\alpha) = \prod_{\upsilon \in S} \|\alpha\|_{\upsilon}$ 
denote the S-norm of $\alpha$.
If $\alpha \in O_S - \{0\}$ then $N_S(\alpha)$ is a 
positive integer and  $N_S(\alpha) = 1$ for $\alpha \in O_S^*$.
For ${\bf x} = (x_0, \dots, x_m) \in {\bf K}^{m+1}$, define the $S$-height as 
$H_S({\bf x}) = \prod_{\upsilon \in S}  \|{\bf x}\|_{\upsilon}$. If 
${\bf x} \in O_S^{m+1} - \{0\}$, then $H_S({\bf x}) \ge 1$ and 
 $H_S(\alpha {\bf x}) = H_S({\bf x})$ for $\alpha \in O_S^*$. 
For a polynomial 
$P$ with coefficients in ${\bf K}$, let $H_S(P)$ denote the $S$-height of 
that point whose coordinates are the coefficients of $P$. 

\bigskip
Linear forms $L_1, \dots, L_q \in {\bar {\bf Q}}[x_0, \dots x_m]$ are said to
 be in {\sl general position} if any $m+1$ of them are linearly independent 
over  ${\bar {\bf Q}}$. 

\bigskip
Let $q, m$ denote positive integers with $q > 2m$, and let
 $F({\bf X}) = F(X_0, \dots, X_m) \in O_S[{\bf X}]$ be a decomposable
 form of degree $q$. For given real numbers $c, \nu$ with $c>0$,
 consider the solutions of the 
inequality 
\begin{equation}
\begin{split}
0 < N_S(F({\bf x})) \leq c H_S({\bf x})^{\nu} ~\mbox{in} ~{\bf x} \in O_S^{m+1}.
\end{split}
\end{equation}
If ${\bf x}$ is a solution of (2.1), then so is 
${\bf x}' = \eta {\bf x}$ for every $\eta \in O^*_S.$ 
Such solution ${\bf x}, {\bf x}'$ are called ${\bf O}^*_S$-proportional.

\bigskip
\noindent{\bf Theorem 1.}
~~{\sl Suppose that $\nu < q-2m$ and that the linear factors of $F$ over ${\bar {\bf Q}}$ are in general position. Then (2.1) has only finitely many 
$O^*_S$-nonproportional solutions. }

\bigskip
\noindent{\bf Remark 1.}
~~Theorem 1 does not remain valid in general for $\nu=q-2m$ (cf. Schmidt 
[Sch 1]). In the important case when ${\bf K} = {\bf Q}$, Theorem 1 
provides an improvement of the previous results of Schmidt [Sch 1, Sch 2] and 
Schlickewei [Schl].  Further, for 
${\bf K} = {\bf Q}, S={\bf M}_{\infty} ({\bf Q})$ and $m=1$, Theorem 1 
gives theorem 3B in Chapter V of Schmidt [Sch 2].

\noindent{\bf Remark 2.}
~For $\nu = 0$, Theorem 1 was proved in a quantitative form in
Gy\H ory [Gy 3]. For a generalization considered over an arbitrary finitely 
generated ground ring over ${\bf Z}$, see  Gy\H ory [Gy 1]. 

\bigskip
The main result of this paper is Theorem 2. Theorem 1 is an immediate 
consequence of this theorem.

\bigskip
\noindent{\bf Theorem 2.}
~~{\sl Let $q, m$ be positive integers with $q > 2m$. Let $c, \nu$ be 
real numbers
 with $c> 0, \nu < q-2m$ and ${\bf G}$ a finite extension of ${\bf K}$.
 For $n=1,2, \dots$, let 
$F_n({\bf X}) = F_n(X_0, \dots, X_m) \in O_S[{\bf X}]$ denote a 
decomposable form of
 degree $q$ which factorizes into linear factors over ${\bf G}$, and suppose 
that these factors are in general position for each $n$.
Then there does not exist an infinite sequence of
 $O^*_S$-nonproportional
 ${\bf x}_n \in O_S^{m+1}, n=1,2, \dots, $ for which 
\begin{equation}
\begin{split}
0 < N_S(F_n({\bf x}_n)) \leq c H_S({\bf x}_n)^{\nu} ~\mbox{for}~ n=1,2, \dots,
\end{split}
\end{equation}
and
\begin{equation}
\begin{split}
 & h(F_n) = o(h({\bf x}_n)) 
~~\mbox{if}~~h({\bf x_n}) \rightarrow \infty~\mbox{as} ~n \rightarrow \infty
\end{split}
\end{equation}
hold.}

\section{Proof of Theorem 2}
We keep the notation of Section 2 and 
recall the Schmidt's subspace theorem with moving targets 
proved by Min Ru and P. Vojta (see [RV]). The fixed target case, 
i.e. the case when $L_{j, n} = L_j$ for each $j$ and $n$,
is due to Min Ru and P. M. Wong (see Theorem 4.1 in [RW]).

\bigskip
\noindent{\bf Theorem A} [RV].
~{\sl Given linear forms
$L_{1,n}, \dots, L_{q,n} \in {\bf K}[X_0, \dots, X_m],
n=1,2, \dots,$ in general 
position for each $n$ and  a sequence 
${\bf x}_n \in {\bf K}^{m+1}$ 
 such that for $j=1,\dots,q$, 
$L_{j,n}({\bf x}_n)\not = 0$ and $h(L_{j,n})= o(h({\bf x}_n))$ 
as $n \rightarrow \infty$.
Then, for any $\epsilon > 0$, 
there exists an infinite subsequence $\{n_k\} \subseteq {\bf N}$ such that
\[ {1\over [{\bf K}:{\bf Q}]} \sum_{\upsilon \in S}\sum_{j=1}^q
\log{ 
\|{\bf x}_{n_k}\|_{\upsilon} \cdot \|L_{j,n_k}\|_{\upsilon} \over 
\|L_{j,n_k}({\bf x}_{n_k})\|_{\upsilon}}    
  \leq (2m  + \epsilon) h({\bf x}_{n_k}) \]
for all $k$.}

\bigskip
The above statement is contained in 
the second part of Theorem 3.1 in [RV]. 
Note that, due to a printing mistake,
 the term $(2m+\epsilon) h({\bf x}_{n_k})$ 
on the right-hand side of 
above inequality 
 was incorrectly stated as 
$(2m+1+\epsilon) h({\bf x}_{n_k})$ in [RV].

\bigskip
\noindent{\sl Proof of Theorem 2.}~
We shall prove Theorem 2 by using the above theorem of [RV].
Assume that there is an infinite sequence ${\bf x}_n = 
(x_{0,n}, \dots, x_{m,n}) \in  O_S^{m+1}$ which  
satisfies (2.2).
First consider the case when the values $h({\bf x}_n)$ are bounded. We may 
assume without loss of generality that $x_{0,n} \not= 0$ for each $n$. Then 
the $h({\bf x}_n /x_{0,n})$ are bounded and this implies that 
${\bf x}_n/x_{0,n}$ may assume only finitely many values in ${\bf K}^{m+1}$. 
Hence there are infinitely many $n$ such that 
${\bf x}_n = x_{0,n} {\bf x}_0$ for some 
${\bf x}_0 \in {\bf K}^{m+1}$. For these $n$ we deduce from (2.2) that 
\[0 <N_S(x_{0,n})^q  N_S(F_n({\bf x}_0)) \leq c N_S(x_{0,n})^{\nu} 
H_S({\bf x}_0)^{\nu}\]
and hence $N_S(x_{0,n})$ are bounded. Since $x_{0,n} \in  O_S$, 
it follows (see e.g. [EGY 3]) that there are infinitely many $n$ for
 which $x_{0,n} = \eta_n x_0'$ with some $\eta_n \in  O^*_S$ and fixed 
$x_0' \in  O_S$. This implies that for these $n$ the ${\bf x}_n$ 
considered above are $O^*_S$-proportional which is a contradiction.

Next consider the case when $h({\bf x}_n)$ are not bounded. We may assume 
that $h({\bf x}_n) \rightarrow \infty$ as $n \rightarrow \infty$. Then, by 
assumption, (2.3) also holds. Further it follows that 
$H_S({\bf x}_n) \rightarrow \infty$ as $n \rightarrow \infty$. 
For $n=1, 2, \dots$, let $F_n = L_{1,n} \dots L_{q,n}$ be a 
factorization of $F_n$ over ${\bf G}$
 into linear factors. Then by Prop. 2.4 of Chapter 3 in [L] it follows that 
\[ \max_j h(L_{j,n}) \leq h(F_n) + c_1\]
where $c_1$ is a positive constant which depends only on $q, m$ and ${\bf G}$.
Together with (2.3) this gives
\begin{equation}
\begin{split}
~\max_jh(L_{j,n}) = o(h({\bf x}_n))~\text{as} ~n \rightarrow \infty.
\end{split}
\end{equation}

Let ${\bf M}({\bf G})$ denote the set of places of ${\bf G}$. For 
$\upsilon 
\in {\bf M}({\bf G})$, define and normalize $\|~\|_{\upsilon}$ in a similar 
manner as over ${\bf K}$ above. Further, let $T$ denote the set of extension 
to ${\bf G}$ of the places in $S$. Then we deduce from (2.2) that 
\begin{equation}
\begin{split}
0 < N_T((F_n({\bf x}_n)) =
N_S(F_n({\bf x}_n))^{[{\bf G}: {\bf K}]} \leq 
(c H_S({\bf x}_n)^{\nu})^{[{\bf G}: {\bf K}]} = c_2 H_T({\bf x}_n)^{\nu},
\end{split}
\end{equation}
where $c_2 = c^{[{\bf G} :{\bf K}]}$. Here $N_T(~), ~H_T(~)$ are defined over 
${\bf G}$ in the same way as $N_S(~), ~H_S(~)$ over ${\bf K}$.

Let $\epsilon > 0$ with $0 < \epsilon < q-2m - \nu$. Then by above 
Theorem A of 
[RV], there is an infinite subsequence 
${\bf x}_{n_k} \in O_S^{m+1}, k=1,2, \dots, $ of 
$\{{\bf x}_n\}$, without loss of generality, we assume $\{{\bf x}_n\}$ itself,
such that 
\[ {1\over [{\bf G}:{\bf Q}]} \sum_{\upsilon \in T}\sum_{j=1}^q 
\log{\|{\bf x}_{n}\|_{\upsilon}\cdot \|L_{j,n}\|_{\upsilon} 
\over \|L_{j,n}({\bf x}_n)\|_{\upsilon} } 
   \leq (2m + \epsilon) h({\bf x}_n). \]
However, $F_n({\bf x}_n) = \prod_{j=1}^q  L_{j, n}({\bf x}_n)$. Furthermore, in view of ${\bf x}_n \in O_S^{m+1}$, we have 
\[h({\bf x}_n)\leq {1\over [{\bf G} : {\bf Q}]} \log H_T({\bf x}_n).\]
Hence it follows that 
\begin{equation}
\begin{split}
&\prod_{\upsilon \in T} 
{\|{\bf x}_{n}\|_{\upsilon}^q\cdot \prod_{j=1}^q \|L_{j,n}\|_{\upsilon} 
\over \|F_n({\bf x}_n)\|_{\upsilon} } 
   \leq H_T({\bf x}_n)^{2m+\epsilon},
\end{split}
\end{equation}
whence
\begin{equation}
\begin{split}
&{H^q_T({\bf x}_n)\cdot \prod_{\upsilon \in T} 
\prod_{j=1}^q \|L_{j,n}\|_{\upsilon} 
\over N_T(F_n({\bf x}_n))} 
   \leq H_T({\bf x}_n)^{2m+\epsilon}. 
\end{split}
\end{equation}
Since the coefficients of $L_{j, n}$ are $T$-integers,  
\begin{equation}
\begin{split}
&\prod_{\upsilon \in T} \prod_{j=1}^q \|L_{j,n}\|_{\upsilon} \ge 1,~ 
\mbox{for}~ n=1,2 \dots \ .  
\end{split}
\end{equation}
Furthermore, it follows from (3.2) that 
\begin{equation}
\begin{split}
& N_T(F_n({\bf x}_n)) \leq c_2 H_T({\bf x}_n)^{\nu}~ \mbox{for}~ n=1,2 \dots \ .
\end{split}
\end{equation}
Combining (3.4), (3.5) and (3.6) gives
\[ H_T({\bf x}_n)^q \leq c_2 H_T({\bf x}_n)^{\nu + 2m + \epsilon}.\]
Since $H_T({\bf x}_n) \rightarrow \infty$ as $n \rightarrow \infty$,
and $q > \nu + 2m + \epsilon$,  this is a
 contradiction.~~~~~~~~~~~~~~~~~~~~~~~~~~~~~~~~~~~~~~~Q.E.D. 

\section{Consequences of Theorem 2}
In this section, we give four further consequences of our Theorem 2.
 
\bigskip
First, it is easy to see that Theorem 2 implies the following result
concerning the S-integer solutions of 
a sequence of decomposable form equations of the form (1.2).

\bigskip
\noindent{\bf Theorem 3.}
~{\sl
 Given positive integers $q, m$ with $q > 2m$, a finite extension 
${\bf G}$ of ${\bf K}$, and 
 a sequence of polynomials $G_n({\bf X}) \in  O_S[{\bf X}]$ 
in ${\bf X} = (X_0, \dots, X_m)$ such that 
 $\deg (G_n({\bf X})) <q-2m$ for $n=1,2, \dots \ .$
Let $F_n({\bf X}) = F_n(X_0,...,X_m) \in  O_S[{\bf X}]$ be 
 a sequence of  decomposable forms of degree q 
such that $F_n$ factorizes into linear forms over ${\bf G}$ which forms are 
in general position for each $n$.
Then there does not exist an infinite sequence of 
$O^*_S$-nonproportional ${\bf x}_n \in O_S^{m+1}$ for which 
\begin{equation}
\begin{split}
F_{n}({\bf x}_{n}) = G_{n}({\bf x}_{n}) \not= 0, ~~n=1,2, \dots,
\end{split}
\end{equation}
\begin{equation}
\begin{split} 
\log H_S(G_n) =  o( \log H_S({\bf x}_{n})), ~\mbox{if} ~H_S({\bf x}_n) 
\rightarrow \infty ~~~\mbox{as}~~ n \rightarrow \infty 
\end{split}
\end{equation} 
and
\begin{equation}
\begin{split}
h(F_{n}) =  o( h({\bf x}_{n})), ~\mbox{if} ~h({\bf x}_n) 
\rightarrow \infty~~as~ n \rightarrow \infty .
\end{split}
\end{equation}}

\bigskip
\noindent{\bf Remark 3.} 
~We note that under the assumptions (4.2) and (4.3), equation (4.1) may have
 infinitely many  $O^*_S$-proportional solutions ${\bf x}_n$. Indeed, 
this is the case if $F_n$ is the same for each $n$, $\eta_n$ is an infinite
 sequence of $S$-units, ${\bf x}_0 \in  O_S^{m+1}$ with
 $F_n({\bf x}_0) = 1, G_n = \eta_n^q$ and 
${\bf x}_n = \eta_n \cdot {\bf x}_0$ for $n=1, 2, \dots \ .$

\bigskip
\noindent{\sl Proof of Theorem 3.}~
Suppose that there is an infinite sequence of 
$O^*_S$-nonproportional ${\bf x}_n \in O^{m+1}_S$ 
satisfying the 
conditions of Theorem 3. 
It is easy to see that the $H_S({\bf x}_n)$ are not bounded. We may assume
 that $H_S({\bf x}_n) \rightarrow \infty$ as 
$n \rightarrow \infty$. 
It follows that 
\begin{equation}
\begin{split}
& 0 < N_S(F_n({\bf x}_{n})) = N_S(G_{n}({\bf x}_{n})) = \prod_{\upsilon \in S}\|G_{n}({\bf x}_{n})\|_{\upsilon} \\
&\leq c_3 \prod_{\upsilon \in S}(\|G_n\|_{\upsilon} \cdot \|{\bf x}_n\|_{\upsilon}^{\deg G_n}) = c_3 H_S(G_n) H_s({\bf x}_n)^{\deg G_n} \\
\end{split}
\end{equation}
where $c_3$ is a positive constant which depends only on $q$ and ${\bf G}$. 
We choose $\nu$ such that $q - 2m - 1 < \nu < q - 2m$. Then, by (4.2), we 
deduce that 
\[ H_S(G_n) \leq H_S({\bf x}_n)^{\nu - \deg G_n}~\text{as}~n \rightarrow \infty.\]
Hence (4.4) implies that 
 \[ 0 < N_S(F_n({\bf x}_n)) \leq c_3 H_S({\bf x}_n)^{\nu}~\text{as}~n \rightarrow \infty\]
which  contradicts Theorem 2. ~~~~~~~~~~~~~~~~~~~~~~~~~~~~~~~~~Q.E.D.

\bigskip
We deduce from Theorem 3 the following. 

\bigskip
\noindent{\bf Theorem 4.}
~{\sl
Given positive integers $q, m$ with $q > 2m$, and 
 a polynomial $G({\bf X}) \in  O_S[{\bf X}]$ 
in ${\bf X}=(X_0, \dots, X_m)$ 
with  total degree less than $q-2m$.
 Let $F({\bf X}) \in  O_S[{\bf X}]$ be a
  decomposable form whose linear factors are in general position. 
Then the equation 
\begin{equation}
\begin{split}
F({\bf X}) = G({\bf X})
\end{split}
\end{equation}
has only finitely many 
solutions
 ${\bf x} = (x_0, \dots, x_m) \in  O_S^{m+1} $ 
with $G({\bf x}) \not= 0$.}

\bigskip
\noindent{\bf Remark 4.}
 For ${\bf K} = {\bf Q}, S={\bf M}_{\infty}({\bf Q})$ and 
$m=1$, our Theorem 4 gives the second assertion of Theorem 3B of
 [Sch 2, Ch. V]. When $G$ is constant, Theorem 4 was proved over 
more general ground rings in K.~Gy\H ory [Gy 1].

\bigskip
\noindent{\sl Proof of Theorem 4.}~
 Suppose that (4.5) has infinitely many solutions 
${\bf x}$ with $G({\bf x}) \not= 0$. Then, by Theorem 3, there are also 
infinitely many solutions ${\bf x}$ such that ${\bf x} = \eta {\bf x}_0$ with 
$\eta \in O^*_S$ and with some fixed ${\bf x}_0 \in  O_S^{m+1}$. 
Then it follows from (4.5) that 
\[ \eta^q F({\bf x}_0) = G(\eta {\bf x}_0) .\]
This can be regarded as an equation of degree $q$ in $\eta$ with leading 
coefficient $F({\bf x}_0) \not = 0$. However, this equation has at most $q$ 
solutions in $\eta$ which proves the assertion. ~~~~~~~~~~~~~~~~~~~~~~~~Q.E.D.

\bigskip
Let $q, m$ be positive integers with $q > 2m$, and 
let $P \in O_S[X]$ be 
a polynomial of degree $q$ without multiple zeros. For given $c>0$ and $\nu,$
consider 
the solutions of the resultant inequality
\begin{equation}
\begin{split}
0 < N_S(\text{Res}(P, Q)) \leq c H_S(Q)^{\nu} ~\text{in} ~{Q}
 \in O_S[X] ~\text{of ~ degree}~ m.
\end{split}
\end{equation}
If $Q$ is a solution then so is $\eta Q$ for each $\eta \in  O^*_S$. 
Such solutions $Q, \eta Q$ are called $O^*_S$-proportional.

\bigskip
\noindent{\bf Theorem 5.}
~{\sl If $\nu < q - 2m$, then (4.6) has only 
 finitely many $O^*_S$-nonproportional solutions.}

\bigskip
\noindent{\bf Remark 5.}
~For  ${\bf K} = {\bf Q}$, Theorem 5 
is an improvement of previous results of E. Wirsing [W], W.M. Schmidt [Sch 1] 
and  H.P. Schlickewei [Schl 1]. Theorem 5 does does not remain valid when 
$\nu=q-2m$ (cf. [Sch 1]). For $\nu = 0$, Theorem 5
was proved in a quantitative form in K. Gy\H ory [Gy 2]. For a generalization 
for polynomials considered over more general ground rings, see K. Gy\H ory [Gy 1].

\bigskip
Theorem 5 is an immediate consequence of the next theorem which will be 
deduced from Theorem 2.

\bigskip
\noindent{\bf Theorem 6.}
~{\sl Let $q, m$ be positive integers with $q > 2m$, $c, \nu$ real numbers 
with $c > 0$, $\nu < q - 2m,$ and ${\bf G}$ a finite extension of ${\bf K}$. 
For every integer $n \ge 1$, let  
 $P_n \in O_S[X]$ denote  a  polynomial of degree $q$ with distinct zeros
in ${\bf G}$.
Then there does not exist a sequence of $O^*_S$-nonproportional 
$Q_n \in O_S[X]$ with $\deg Q_n = m$ for which 
\begin{equation}
\begin{split}
& 0 < N_S(\text{Res}(P_n, Q_n)) \leq c H_S(Q_n)^{\nu}, \ n=1,2, \dots , 
\end{split}
\end{equation}
and
\begin{equation}
\begin{split}
&h(P_n) = o(h(Q_n))~~\mbox{if}~h(Q_n) \rightarrow \infty ~\mbox{as} ~n \rightarrow \infty
\end{split}
\end{equation}
hold.}

\bigskip
This should be compared with Corollary 4 of K. Gy\H ory [Gy 2] where 
$\nu = 0$, but $O_S$ is replaced by a more general ground ring.

\bigskip
\noindent{\sl Proof of Theorem 6.}~
Put 
$$
P_n(X) = a_{0,n}(X-\alpha_{1,n}) \dots (X-\alpha_{q,n}) ~ \text{for}~ n=1,2, \dots \ .
$$
Assume that there is   an infinite sequence of 
$$
Q_n(X) = x_{0,n}X^m + x_{1,n}X^{m-1} + \cdots + x_{m,n},
$$
satisfying (4.7) and (4.8).
For $n \geq 1$ set
\[ F_n({\bf X}) = F_n(X_0, \dots, X_m) = a_{0,n}^m \prod_{i=1}^q (X_0 \alpha_{i, n}^m + X_1\alpha_{i,n}^{m-1} + \dots + X_m). \]
Then $F_n$ has its coefficients in $O_S$.
Further it factorizes into linear factors over ${\bf G}$ and these linear
 forms are in general position for every $n$. 
For ${\bf x}_n = (x_{0,n}, \dots, x_{m,n})$ we have 
\[ \text{Res}(P_n, Q_n) = F_n({\bf x}_{n}). \]
 Hence (4.7) 
implies that 
\begin{equation}
\begin{split}
 0 < N_S(F_n({\bf x}_n)) \leq c H_S(Q_n)^{\nu}  =  c H_S({\bf x}_n)^{\nu}.
\end{split}
\end{equation}
Further, by using Prop. 2.4 of Ch. 3 in [L] it is easy to see that 
\begin{equation}
\begin{split}
h(F_n) \leq q h(P_n) + c_4,
\end{split}
\end{equation}
where $c_4$ is a constant which depends only on $q, m$ and ${\bf G}$. Hence if
$h(Q_n) \rightarrow \infty$ as $n \rightarrow \infty$,
it follows from 
(4.10) and (4.8) that
$$ 
h(F_n)= o(h({\bf x}_n)), ~\text{if}~ h({\bf x}_n) \rightarrow \infty~
 \text{as}~ n \rightarrow \infty.
$$ 
Together with (4.9) this contradicts our Theorem 2.~~~~~~~~~~~~~~~~~~~~~~~~~~~
Q.E.D.

\section{Some generalizations}
In the previous sections, we assume that $q > 2m$ where $q$ is the degree of 
a decomposable form $F$, and $m+1$ is the number of variables in $F$. In this 
section, we consider the case $q > 2m -l +1$, where $l$ is an integer with 
$ 1 \leq l \leq m+1$.
As we indicated earlier, in the case that $l > 1$, 
finiteness result is not expected. Rather,
we show that the set of solutions is contained in a finite union of 
proper subspaces in this situation.

\bigskip
First of all, Theorem 1 can be generalized as follows.

\bigskip
\noindent{\bf Theorem 7.}
{\sl Let $q, m, l$ denote positive integers with 
$q > 2m - l + 1, 1 \leq l \leq m+1$, and 
${\bf G}$ a finite extension of ${\bf K}$.  Let
 $F({\bf X}) = F(X_0, \dots, X_m) \in O_S[{\bf X}]$ be a decomposable
 form of degree $q$ which factorizes into linear factors over ${\bf G}$. 
For given real numbers $c, \nu$ with $c>0$,
 consider the solutions of the 
inequality 
\begin{equation}
\begin{split}
0 < N_S(F({\bf x})) \leq c H_S({\bf x})^{\nu} ~\mbox{in} ~{\bf x} \in O_S^{m+1}.
\end{split}
\end{equation}
Suppose that $\nu < q - 2m + l - 1$ and that the linear
factors of $F$ over ${\bar {\bf Q}}$ are in general position. Then the set of
solutions of  (5.1) is contained in a finite union of linear subspaces of 
${\bf K}^{m+1}$ of dimension at most $l$.}

It is easy to see that for $l=1$, Theorem 7 
gives Theorem 1. For $\nu=0$, this can be compared with Theorem 3 of [EGy1]

\bigskip
The proof of Theorem 7 is similar to the proof of Theorem 2. First we recall
Theorem 4.1 of Ru-Wong [RW] (see also [RV], Theorem 3.1).

\bigskip
\noindent{\bf Theorem B} (Theorem 4.1 of [RW]).
~{\sl Given linear forms
$L_1, \dots, L_q \in {\bf K}[X_0, \dots, X_m]$ in general 
position.  Then for any $\epsilon > 0$, the set of points 
${\bf x} \in  {\bf K}^{m+1}$ 
such that $L_j({\bf x}) \ne 0$ for $j=1, \dots ,q$ and
\[ {1\over [{\bf K}:{\bf Q}]} \sum_{v \in S}\sum_{j=1}^q
\log{ 
\|{\bf x}\|_{\upsilon} \cdot \|L_j\|_{\upsilon} \over 
\|L_j({\bf x})\|_{\upsilon}}    
  \ge (2m - l + 1 + \epsilon) h({\bf x}) \]
 is contained in a finite union of linear subspaces of ${\bf K}^{m+1}$ 
of dimension at most $l$.}

\bigskip
\noindent{\sl Proof of Theorem 7.}~
Let ${\bf x} \not= {\bf 0}$ be a solution of (5.1). 
Let $F = L_1 \dots L_q$ be a 
factorization of $F$ over ${\bf G}$
 into linear factors. 
Let ${\bf M}({\bf G})$ denote the set of places of ${\bf G}$. For 
$\upsilon 
\in {\bf M}({\bf G})$, define and normalize $\|~\|_{\upsilon}$ in a similar 
manner as over ${\bf K}$ above. Further, let $T$ denote the set of extension
to ${\bf G}$ of the places in $S$. Then we deduce from (5.1) that 
\[ 0 < N_T((F({\bf x})) =
N_S(F({\bf x}))^{[{\bf G}: {\bf K}]} \leq 
(c H_S({\bf x})^{\nu})^{[{\bf G}: {\bf K}]} = c_2 H_T({\bf x})^{\nu},\]
where $c_2 = c^{[{\bf G} :{\bf K}]}$. Here $N_T(~), ~H_T(~)$ are defined over 
${\bf G}$ in the same way as $N_S(~), ~H_S(~)$ over ${\bf K}$.
We have $F({\bf x}) = \prod_{j=1}^q  L_{j}({\bf x})$.
 It follows that 
\begin{equation}
\begin{split}
&\prod_{\upsilon \in T} 
{\|{\bf x}\|_{\upsilon}^q\cdot \prod_{j=1}^q \|L_j\|_{\upsilon} 
\over \prod_{j=1}^q  \|L_j({\bf x})\|_{\upsilon} }  = {
H_T({\bf x})^q \cdot 
(\prod_{\upsilon \in T} \prod_{j=1}^q \|L_j\|_{\upsilon}) 
\over N_T(F({\bf x})) } \\
&\ge {H_T({\bf x})^q \cdot 
(\prod_{\upsilon \in T} \prod_{j=1}^q \|L_j\|_{\upsilon}) 
\over c_2 H_T({\bf x})^{\nu}} 
= c_5 H_T({\bf x})^{q - \nu},\\
\end{split}
\end{equation}
where $c_5 = \prod_{\upsilon \in T}  \prod_{j=1}^q \|L_j\|_{\upsilon}/c_2$.
In view of ${\bf x} \in O^m_T$, we have 
$$
h({\bf x}) \leq {1 \over [{\bf G}:{\bf Q}]} \log H_T({\bf x}).
$$
Thus, by taking logarithms on both side of (5.2), 
\[{1\over [{\bf G} : {\bf Q}]} \sum_{\upsilon \in T} \sum_{j=1}^q  
\log {\|{\bf x}\|_{\upsilon}\cdot  \|L_j\|_{\upsilon} 
\over  \|L_{j}({\bf x})\|_{\upsilon} } 
\ge (q - \nu)h({\bf x}) + c_6. \]
where $c_6=\log c_5 / [{\bf G}:{\bf Q}]$. 
Let $\epsilon > 0$ with $ q - \nu \geq 2m + l - 1 > 2 \epsilon$, then  
\[{1\over [{\bf G} : {\bf Q}]} \sum_{\upsilon \in T} \sum_{j=1}^q 
\log {\|{\bf x}\|_{\upsilon}\cdot  \|L_j\|_{\upsilon} 
\over  \|L_{j}({\bf x})\|_{\upsilon} } 
\ge (2m -l + 1 + 2\epsilon) h({\bf x}) + c_6. \]
Since the set of points ${\bf x} \in {\bf G}^{m+1}$ with
 $\epsilon h({\bf x}) + c_6 < 0$ is finite, 
excluding these points yields  
\[{1\over [{\bf G} : {\bf Q}]} \sum_{\upsilon \in T} \sum_{j=1}^q 
\log {\|{\bf x}\|_{\upsilon}\cdot  \|L_j\|_{\upsilon} 
\over   \|L_{j}({\bf x})\|_{\upsilon}} 
\ge (2m - l + 1 + \epsilon)h({\bf x}). \]
Thus, the above quoted Theorem B of [RW] implies Theorem 7.  ~~~~~~~~~~~~~~~~~~~~~~~~~~~~~~~~~~~Q.E.D.

\bigskip
Similarly, by using Theorem 3.1 of [RV] Theorem 2 can be generalized
as follows:

\bigskip
\noindent{\bf Theorem 8.}
~~{\sl Let $q, m, l$ be positive integers with $q > 2m -l + 1, 
1 \leq l \leq m+1$. 
Let $c, \nu$ be 
real numbers
 with $c> 0, \nu < q-2m + l -1$ and ${\bf G}$ a finite extension of ${\bf K}$.
 For $n=1,2, \dots$, let 
$F_n({\bf X}) = F_n(X_0, \dots, X_m) \in O_S[{\bf X}]$ denote a 
decomposable form of
 degree $q$ which factorizes into linear factors over ${\bf G}$, and suppose 
that these factors are in general position for each $n$.
Then the points
${\bf x}_n \in O_S^{m+1}, n=1,2, \dots, $ satisfying
\[
0 < N_S(F_n({\bf x}_n)) \leq c 
H_S({\bf x}_n)^{\nu} \ \text{for} \  n=1,2, \dots, \]
and
\[
 h(F_n) = o(h({\bf x}_n)) 
~~\mbox{if}~~h({\bf x_n}) \rightarrow \infty~\mbox{as} ~n \rightarrow \infty,\]
are $l$-degenerate. For the concept of $l$-degenerate, see [RV].}

Theorems 7 and 8 have applications of the same type to resultant 
inequalities as Theorems 1 and 2 above.

\bigskip
As was shown earlier, Theorem 1 implies Theorem 4. Similarly,
our Theorem 7 has the 
following consequence. 

\bigskip
\noindent{\bf Theorem 9.}
~{\sl
Given positive integers $q, m, l$ with $q > 2m - l + 1$, $1 \leq l \leq m+1$,
 and 
 a polynomial $G({\bf X}) \in  O_S[{\bf X}]$ 
in ${\bf X}=(X_0, \dots, X_m)$ 
with  total degree less than $q-2m + l -1$.
 Let $F({\bf X}) \in O_S[{\bf X}]$ be a
decomposable form of degree q whose linear factors are in general position. 
Then the set of solutions of
\[
F({\bf x}) = G({\bf x})~~~~\mbox{in}~{\bf x } = 
(x_0, \dots, x_m) \in O_S^{m+1} \]
is contained in a finite union of linear subspaces of ${\bf K}^{m+1}$ of 
dimension at most $l$.}

\bigskip

\section{References} 
\vspace{4mm}
\noindent
[EGy 1] J.H. Evertse and K. Gy\H ory, {\sl Finiteness criteria for decomposable form equations}, 
Acta Arith., 50(1988), 357-379. 

\vspace{4mm}
\noindent
[EGy 2] J.H. Evertse and K. Gy\H ory, {\sl Decomposable form equations,} 
in ``New Advances in Transcendence Theory'' (A. Baker ed.), pp. 175-202. 
Cambridge University Press, 1988.

\vspace{4mm}
\noindent
[EGy 3] J.H. Evertse and K. Gy\H ory, {\sl Effective finiteness results for 
binary forms with given discriminant}, Compositio Math., 79(1991), 169-204.

\vspace{4mm}
\noindent
[EGy 4] J.H. Evertse and K. Gy\H ory, {\sl The number of families of solutions of decomposable form equations}, Acta Arith., 80(1997), 367-394. 

\vspace{4mm}
\noindent
[Gy 1] K. Gy\H ory, {\sl Some applications of decomposable form equations to 
resultant equations}, 
Colloq. Math., 65(1993), 267-275.

\vspace{4mm}
\noindent
[Gy 2] K. Gy\H ory, {\sl On the number of pairs of polynomials with given 
resultant or given semi-resultant}, Acta Sci. Math., 57(1993), 515-529.

\vspace{4mm}
\noindent
[Gy 3] K. Gy\H ory, {\sl On the irreducibility of neighbouring polynomials},
 Acta Arith., 67(1994), 283-296.

\vspace{4mm}
\noindent
[L] S. Lang, {\sl Fundamentals of Diophantine Geometry}, Springer Verlag, 
1983.

\vspace{4mm}
\noindent
[RV] Min Ru and P. Vojta, 
{\sl Schmidt's subspace theorem with moving targets},
 Invent. Math., 127(1997), 51-65.

\vspace{4mm}
\noindent
[RW] Min Ru and P. M. Wong, {\sl Integral points of
${\bf P}^n - \{2n+1$ hyperplanes in general position$\}$},
Invent. Math., 106 (1991), 195-216.

\vspace{4mm}
\noindent
[Schl] H.P. Schlickewei, {\sl Inequalities for decomposable forms}, Ast\'erisque 41-42(1977), 267-271.

\vspace{4mm}
\noindent
[Sch 1] W. M. Schmidt,
{\sl Inequalities for resultants and for decomposable forms}, In:
{\sl Diophantine approximation and its applications},
 Academic. Press, New York, 1973, 235-253.

\vspace{4mm}
\noindent
[Sch 2] W. M. Schmidt,
{\sl Diophantine Approximation}, Lecture Notes Math.\ 785,
Springer Verlag, Berlin etc., 1980.

\vspace{4mm}
\noindent
[W] E. Wirsing, {\sl On approximations of algebraic numbers by algebraic 
numbers of bounded degree}, Proc. of Symp. in Pure Math., XX(1971), 
pp. 213-247.

\end{document}